\documentclass[12pt, reqno]{amsart}
\usepackage{amsmath, amsthm, amscd, amsfonts, amssymb, graphicx, color,verbatim}

\newtheorem{thm}{Theorem}[section]
\newtheorem{cor}[thm]{Corollary}

\newtheorem{prop}[thm]{Proposition}
\theoremstyle{definition}

\theoremstyle{remark}
\newtheorem{rem}[thm]{Remark}

\newtheorem{ex}[thm]{Example}
\newcommand{\bpf}{\begin{proof}}
\newcommand{\epf}{\end{proof}}

\numberwithin{equation}{section}
\newcommand{\C}{{\mathbb C}}
\newcommand{\B}{{\mathbb B}}
\newcommand{\D}{{\mathbb D}}
\newcommand{\T}{{\mathbb T}}

\newcommand{\clh}{\mathcal{H}}

\newcommand{\clm}{\mathcal{M}}

\newcommand{\raro}{\rightarrow}

\newcommand{\vp}{{\varphi}}

\begin{document}

\setcounter{page}{1}

\today

\title[New norm estimate for composition operators]{New norm estimate for composition operators}

\author[Preeti, Muthukumar \MakeLowercase{and} Sarkar]{Preeti Kumari, P. Muthukumar \MakeLowercase{and} Jaydeb Sarkar}

\address{Preeti Kumari, Department of Mathematics and Statistics, Indian Institute of Technology,
 Kanpur- 208016, India.}
\email{preetisharma8274@gmail.com, preetik22@iitk.ac.in}

\address{P. Muthukumar, Department of Mathematics and Statistics, Indian Institute of Technology,
 Kanpur- 208016, India.}
\email{pmuthumaths@gmail.com, muthu@iitk.ac.in}

\address{Jaydeb Sarkar, Indian Statistical Institute,
Statistics and Mathematics Unit \\
 8th Mile, Mysore Road\\
Bangalore 560 059\\ India.}
 \email{jaydeb@gmail.com, jay@isibang.ac.in}

\subjclass[2020]{47B33, 30H10, 42C99, 46E22}

\keywords{Composition operators, Hardy space, norm, reproducing kernels, bounded analytic functions}

\begin{abstract}
The classical Littlewood's theorem establishes boundedness and provides a norm estimate for composition operators on the Hardy space. In this paper, we offer an alternative proof of boundedness and derive a new norm estimate that improves upon the classical bound given by Littlewood’s theorem.
\end{abstract}

\maketitle

\tableofcontents

\section{Introduction}\label{sec: intro}

We denote the open unit disc in the complex plane by $\D$, and let $H^2(\D)$ denote the Hardy space over $\D$ \cite{shap}. For a self-map $\vp$ of $\D$ (where self-maps are always assumed to be analytic), the \textit{composition operator} $C_\vp$ is defined by
\[
C_\vp f = f \circ \vp,
\]
for all $f \in H^2(\D)$. The classical Littlewood's theorem (see Shapiro, \cite[page 16]{shap}) asserts that $C_\vp$ is a bounded linear operator on $H^2(\D)$. Moreover, its operator norm satisfies the estimate
\[
\|C_\vp\| \leq \sqrt{\frac{1 + |\vp(0)|}{1 - |\vp(0)|}}.
\]
This result is based on Littlewood's celebrated and century-old subordination principle \cite{lit}, which was further fine tuned to bring up the full form as described above by Nordgren \cite{Nord} and Ryff \cite{ryff} in the late 1960s.

The purpose of this paper is to improve the above classical norm estimate. More specifically, we prove (see  Theorem \ref{Bdd}):

\begin{thm}
Let $\vp$ be a self-map of $\D$. Then
\[
\|C_{\vp}\|\leq \left\|{\frac{\sqrt{1-|\vp(0)|^2}}{1-\overline{\vp(0)} \vp}}\right\|_\infty \leq \sqrt{\frac{1+|\vp(0)|}{1-|\vp(0)|}}.
\]
\end{thm}

Along the way, we also present an alternative proof of Littlewood’s theorem establishing the boundedness of composition operators on $H^2(\D)$. In this paper, we follow standard notation: for an analytic function $\psi: \D \raro \mathbb{C}$, we write
\[
\|\psi\|_\infty = \sup_{z \in \D} |\psi(z)|,
\]
and denote the space of all such $\psi$ with $\|\psi\|_\infty < \infty$ by $H^\infty(\D)$.

We prove, moreover, that the new estimate is sharper for a broader class of self-maps. For instance, if
\[
\vp(0) \neq 0,
\]
and
\[
\|\vp\|_\infty < 1,
\]
then
\[
\left\|{\frac{\sqrt{1-|\vp(0)|^2}}{1-\overline{\vp(0)}\vp}}\right\|_\infty < \sqrt{\frac{1+|\vp(0)|}{1-|\vp(0)|}}.
\]
We remark that self-maps $\varphi$ satisfying $\|\varphi\|_\infty < 1$ induce compact composition operators. From this perspective, we also provide examples of self-maps $\varphi$ with $\|\varphi\|_\infty = 1$ for which the new norm estimate remains sharper than the classical one. Notably, such maps may also give rise to noncompact composition operators. Furthermore, we classify all self-maps $\varphi$ for which the new and the classical estimates coincide. In fact, this occurs, in particular, when $\varphi(0) = 0$, in which case the classical estimate, the new estimate, and the norm of $C_\varphi$ all coincide and are equal to $1$. The remaining part is presented in the following result (see Theorem \ref{equal}).

\begin{thm}
Let $\vp$ be a self-map of $\D$. Assume that $\vp(0) \neq 0$. Then
\[
\left\|{\frac{\sqrt{1-|\vp(0)|^2}}{1-\overline{\vp(0)}\vp}}\right\|_\infty = \sqrt{\frac{1+|\vp(0)|}{1-|\vp(0)|}},
\]
if and only if
\[
\dfrac{\vp(0)}{|\vp(0)|}\in \partial( \vp(\D)).
\]
\end{thm}

The central idea behind the new estimate presented here lies in the application of reproducing kernel Hilbert space techniques. This framework naturally extends to more general contexts, allowing us to establish boundedness and obtain norm estimates for composition operators on reproducing kernel Hilbert spaces. The general result is presented in Theorem \ref{genbdd}.

We also address a statement incorrectly asserted in Schwartz’s 1969 thesis. In particular, Schwartz claimed that if $\vp$ is a self-map of $\D$ with $\|\vp\|_\infty<1$ and $\vp(0)\neq0$, then
\[
\|C_{\vp}\| < \sqrt{\dfrac{\|\vp\|_\infty+|\vp(0)|}{\|\vp\|_\infty-|\vp(0)|}}.
\]
He further asserted that this estimate for the norm of $C_\vp$ is sharper than the classical estimate given by Littlewood’s theorem. However, in view of both his and our estimates along with Proposition \ref{prop: Schartz}, we obtain the following chain of inequalities:
\[
\|C_{\vp}\|\leq \left\|{\frac{\sqrt{1-|\vp(0)|^2}}{1-\overline{\vp(0)} \vp}}\right\|_\infty < \sqrt{\frac{1+|\vp(0)|}{1-|\vp(0)|}} < \sqrt{\dfrac{\|\vp\|_\infty+|\vp(0)|}{\|\vp\|_\infty-|\vp(0)|}},
\]
for all self-maps $\vp$ with $\|\vp\|_\infty < 1$ and $\vp(0) = 0$. The above chain of inequalities shows that Schwartz’s assertion that his norm estimate for $C_\vp$ improves upon the classical bound of Littlewood’s theorem is, in fact, incorrect.

The remainder of the paper is organized as follows. Section \ref{sec: CO on H} presents the new norm estimate for composition operators and establishes their boundedness. Section \ref{sec: estimate} proves that the new estimate is sharper than the classical one and characterizes the self-maps for which the two estimates coincide. Section \ref{sec-RKHS} extends the discussion to composition operators in the framework of reproducing kernel Hilbert spaces. Section \ref{sec: examples} provides examples of self-maps illustrating the sharper nature of the new norm estimate. Finally, in Section \ref{sec: Schwartz}, we examine Schwartz’s assertion on the norm estimates for certain self-maps of $\D$.

\section{Composition operators on $H^2(\D)$}\label{sec: CO on H}

The purpose of this section is to prove the central result of this paper: the boundedness of composition operators together with a new norm estimate. The proof relies on three elementary facts about reproducing kernel Hilbert spaces. Along with recalling these three facts, we briefly review the essential reproducing kernel techniques needed for this proof and refer the reader to the monograph \cite{paul} for a more detailed discussion.

Let $X$ be a nonempty set. A function $k: X \times X \rightarrow \mathbb{C}$ is said to be a \textit{kernel function} on $X$ if
\[
\sum\limits_{i,j=1}^m \overline{c_i} c_j k(x_i,x_j) \geq 0,
\]
for all choices of scalars $c_1, \ldots, c_m \in \mathbb{C}$, points $x_1, \ldots, x_m \in X$, and integers $m \in \mathbb{N}$. Such a kernel function is often denoted symbolically by
\[
k(x,y) \geq 0.
\]
Given a kernel function $k$ on a set $X$, we denote by $\clh(k)$ the reproducing kernel Hilbert space corresponding to $k$. For example, the Hardy space is a reproducing kernel Hilbert space, and its kernel function is the Szeg\"{o} kernel $\mathbb{S}$ of the unit disc:
\[
\mathbb{S}(z, w) = \frac{1}{1 - z \bar{w}} \qquad (z, w \in \D).
\]
For each $y \in X$, define a function $k(\cdot, y): X \rightarrow \C$ by
\[
(k(\cdot, y))(x) = k(x, y),
\]
for all $x \in X$. It is easy to see that the set $\{k(\cdot, y): y \in X\}$ is total in $\clh(k)$.

The characterization of the elements of $\clh(k)$ is one of the keys to our methodology. Let $f: X \rightarrow \C$ be a function. Then $f\in \clh(k)$ with
\[
\|f\| \leq c,
\]
if and only if (cf. \cite[Theorem 3.11]{paul})
\begin{equation}\label{eqn: f in RKHS}
c^2 k(x,y) - f(x)\overline{f(y)} \geq 0.
\end{equation}
Let $\tilde{k}$ be another kernel function on $X$. A function $\vp : X \raro \mathbb{C}$ is said to be a \textit{multiplier} from $\clh(k)$ into $\clh(\tilde k)$ if
\[
\vp f \in \clh(\tilde{k}),
\]
for all $f \in \clh(k)$, equivalently,
\[
\vp \clh(k) \subseteq \clh(\tilde{k}).
\]
The set of all such multipliers is denoted by $\clm(\clh(k), \clh(\tilde{k}))$. Note that for a multiplier $\vp \in \clm(\clh(k), \clh(\tilde{k}))$, the \textit{multiplication operator} $M_\vp$ defines a bounded linear operator from $\clh(k)$ to $\clh(\tilde{k})$, where
\[
(M_\vp f)(x) = \vp(x) f(x),
\]
for all $f \in \clh(k)$ and $x \in X$. We need the following characterization of multipliers, which will serve as the second key ingredient in what follows: Suppose $\vp: X \raro \mathbb{C}$ be a function. Then $\vp \in \clm(\clh(k), \clh(\tilde{k}))$ if and only if there exists a constant $c \geq 0$ such that (cf. \cite[Theorem 5.21]{paul})
\begin{equation}\label{eqn: mult}
c^2 \tilde{k}(x, y) - \vp(x) k(x,y) \overline{\vp(y)} \geq 0.
\end{equation}
Moreover, if $\vp \in \clm(\clh(k), \clh(\tilde{k}))$, then
\begin{equation}\label{eqn: norm Mphi}
\|M_\vp\| = \inf \{c \geq 0: c^2 \tilde{k}(x, y) - \vp(x) k(x,y) \overline{\vp(y)} \geq 0\}.
\end{equation}
Finally, we turn to the third and final key point needed for the proof that follows: namely, that $k_\psi(z, w) \ge 0$ on $\mathbb{D}$ for every $\psi \in H^\infty(\mathbb{D})$ with $\|\psi\|_\infty \leq 1$, where $k_\psi$ is defined below. Recall that $\clm(H^2(\D), H^2(\D)) = H^\infty(\D)$. Let $\psi \in H^\infty(\D)$. We denote by $T_\psi = M_\psi$ the analytic Toeplitz operator on $H^2(\D)$ with symbol $\psi$. We know that
\[
\|T_\psi\| = \|\psi\|_\infty.
\]
In particular, if $\|\psi\|_\infty \leq 1$, then $T_\psi$ is a contraction, and then \eqref{eqn: norm Mphi} implies
\[
k_\psi(z, w):= \mathbb{S}(z, w)-\psi(z)\mathbb{S}(z, w)\overline{\psi(w)}=\frac{1- \psi(z) \overline{\psi(w)}}{1-z \overline{w}}\geq 0,
\]
that is,
\[
k_\psi(z, w)= \frac{1- \psi(z) \overline{\psi(w)}}{1-z \overline{w}}\geq 0.
\]
Note that
\[
k_\psi(z, w)= (1- \psi(z) \overline{\psi(w)}) \mathbb{S}(z,w).
\]
The reproducing kernel Hilbert space $\clh(k_\psi)$ corresponding to the kernel $k_\psi$ is the well-known de Branges–Rovnyak space \cite{Sarason}. In the following, we will not use any tools from de Branges–Rovnyak spaces, but only the basic properties of the kernel function $k_\psi$.

\begin{thm}\label{Bdd}
Let $\vp$ be a self-map of $\D$. Then $C_\vp$ is a bounded linear operator on $H^2(\D)$, and
\[
\|C_{\vp}\|\leq \left\|{\frac{\sqrt{1-|\vp(0)|^2}}{1-\overline{\vp(0)}\vp}}\right\|_\infty \leq \sqrt{\frac{1+|\vp(0)|}{1-|\vp(0)|}}.
\]
\end{thm}
\bpf
As we have just observed, $k_\vp$ is a kernel function on $\D$, where
\[
k_\vp(z,w)= \frac{1-\vp(z) \overline{\vp(w)}}{1- z \overline{w}},
\]
for all $z, w \in \D$. The reproducing kernel Hilbert space corresponding to $k_\vp$ is denoted by $\clh(k_\vp)$. Note that
\[
k_\vp(\cdot, 0) = 1- \overline{\vp(0)} \vp \in \clh(k_\vp).
\]
By \eqref{eqn: f in RKHS}, we have $\|k_\vp(\cdot, 0)\|^2_{\clh(k_\vp)} k_\vp(z,w) - k_\vp(z, 0) \overline{k_\vp(w , 0)} \geq0$, that is,
\begin{equation}\label{R1}
\|k_\vp(\cdot, 0)\|^2_{\clh(k_\vp)} \frac{1-\vp(z)\overline{\vp(w)}}{1-z \overline{w}} - k_\vp(z, 0) \overline{k_\vp(w , 0)} \geq0.
\end{equation}
On the other hand, for each $z \in \D$, $k_\vp(z, 0) + \overline{\vp(0)} \vp(z) = 1$ implies
\[
1 \leq |k_\vp(z, 0)| + |\overline{\vp(0)}| |\vp(z)| \leq |k_\vp(z, 0)| + |\overline{\vp(0)}|,
\]
as $\vp$ is a self-map of $\D$. In other words, there exists $\delta > 0$ such that
\[
|k_\vp(z, 0)| \geq 1 - |\vp(0)| > \delta,
\]
for all $z \in \D$ (as $|\vp(0)| < 1$), and consequently
\[
\frac{1}{k_\vp(\cdot, 0)} = \frac{1}{1-\overline{\vp(0)}\vp} \in H^\infty(\D).
\]
In particular, we have
\[
\frac{1}{k_\vp(z, 0) \overline{k_\vp(w, 0)}} \geq 0.
\]
Moreover, since $|\vp(z)| < 1$ for all $z \in \D$, we have the kernel function
\[
\frac{1}{1 - \vp(z) \overline{\vp(w)}} \geq 0,
\]
and hence
\[
\frac{1}{k_\vp(z, 0)} \; \frac{1}{1 - \vp(z) \overline{\vp(w)}} \; \frac{1}{\overline{k_\vp(w, 0)}} \geq 0.
\]
In the above, we have used the fact that the product of two kernel functions is again a kernel function, which follows from the Schur product theorem.
Then, by multiplying the kernel function defined above with the kernel function given in \eqref{R1}, we find that
\begin{equation}\label{R2}
\|k_\vp(\cdot, 0)\|^2_{\clh(k_\vp)} \frac{1}{k_\vp(z, 0)} \; \frac{1}{1 - z \overline{w}} \; \frac{1}{\overline{k_\vp(w, 0)}} - \frac{1}{1 - \vp(z) \overline{\vp(w)}}\geq0.
\end{equation}
Next, we treat $\frac{1}{k_\vp(\cdot, 0)} \in H^\infty(\D)$ as a multiplier of $H^2(\D)$, and find that this leads to (see the classification of multipliers in \eqref{eqn: mult})
\[
\left\|\frac{1}{k_\vp(\cdot, 0)}\right\|_{\infty}^2 \frac{1}{1- z \overline{w}} - \frac{1}{k_\vp(z, 0)} \; \frac{1}{1 - z \overline{w}} \; \frac{1}{\overline{k_\vp(w, 0)}} \geq 0.
\]
First, we multiply this by the positive constant $\|k_\vp(\cdot, 0)\|^2_{\clh(k_\vp)}$, and then, by adding it to \eqref{R2}, we obtain
\[
c^2 \frac{1}{1 - z \overline{w}} - \frac{1}{1 - \vp(z) \overline{\vp(w)}}\geq 0,
\]
where $c > 0$ is given by
\[
c = \|k_\vp(\cdot, 0)\|_{\clh(k_\vp)} \left\|\frac{1}{k_\vp(\cdot, 0)}\right\|_{\infty}.
\]
We claim that the above kernel property implies that $C_\vp$ is bounded on $H^2(\D)$ with $\|C_\vp\| \leq c$. Although this follows directly from \cite[ page 72, Theorem 5.10]{paul}, we present a self-contained proof below: For $f \in H^2(\D)$, we know from \eqref{eqn: f in RKHS} that
\[
\|f\|^2\frac{1}{1 - z \overline{w}} - f(z) \overline{f(w)} \geq 0.
\]
In particular, we have
\[
f(\vp(z)) \overline{f(\vp(w))} \leq \|f\|^2 \frac{1}{1 - \vp(z) \overline{\vp(w)}},
\]
for all $z, w \in \D$. Together with the fact that $c^2 \frac{1}{1 - z \overline{w}} - \frac{1}{1 - \vp(z) \overline{\vp(w)}}\geq 0$, this implies that
\[
f(\vp(z)) \overline{f(\vp(w))} \leq c^2 \|f\|^2 \frac{1}{1 - z \overline{w}},
\]
that is,
\[
c^2 \|f\|^2 \frac{1}{1 - z \overline{w}} - f(\vp(z)) \overline{f(\vp(w))} \geq 0.
\]
In this case, once again, \eqref{eqn: f in RKHS} implies that $f \circ \vp \in H^2(\D)$ and
\[
\|f \circ \vp\| \leq c \|f\|,
\]
which completes the proof of the claim that $C_\vp$ is bounded on $H^2(\D)$ with $\|C_\vp\|\leq c$. Finally, we have
\[
\|k_\vp(\cdot, 0)\|_{\clh(k_\vp)}^2 = k_\vp(0, 0) = 1 - |\vp(0)|^2,
\]
and hence
\[
c = \|k_\vp(\cdot, 0)\|_{\clh(k_\vp)} \left\|\frac{1}{k_\vp(\cdot, 0)}\right\|_{\infty} = \left\|{\frac{\sqrt{1-|\vp(0)|^2}}{1-\overline{\vp(0)}\vp}}\right\|_\infty.
\]
As $\vp$ is a self-map of $\D$, we have
\begin{equation}\label{eqn: 1 - phi bar phi}
|1-\overline{\vp(0)}\vp(z)| \geq 1-|\overline{\vp(0)}|=1-|\vp(0)|,
\end{equation}
for all $z \in \D$. This immediately implies
\[
\left\|{\frac{\sqrt{1-|\vp(0)|^2}}{1-\overline{\vp(0)}\vp}}\right\|_\infty\leq \sqrt{\frac{1+|\vp(0)|}{1-|\vp(0)|}},
\]
and concludes the proof of the theorem.
\epf

The use of reproducing kernel techniques in the theory of composition operators is common. In fact, in \cite{jury}, Jury used kernel function techniques to provide an alternative proof of the boundedness of composition operators. His method, however, relies on a more detailed structure of kernel functions and, in particular, makes essential use of the de Branges–Rovnyak space techniques \cite{Sarason}. In the present situation, we also make contact with the de Branges–Rovnyak kernel function (the kernel function $k_\vp$ used in the proof of Theorem \ref{Bdd}), but we use it only in its capacity as a kernel function, without appealing to the full de Branges–Rovnyak space framework. The present method is not only more direct than that of \cite{jury} but also yields new information concerning the norm estimate.

The computation of norms of composition operators is generally a case dependent problem. In other words, the best one can usually hope for is to obtain general bounds for the norms of all composition operators. The classical estimate
\[
\|C_\vp\| \leq \sqrt{\frac{1+|\vp(0)|}{1-|\vp(0)|}},
\]
is well known. The new estimate obtained above, namely,
\[
\|C_\vp\| \leq \left\|{\frac{\sqrt{1-|\vp(0)|^2}}{1-\overline{\vp(0)}\vp}}\right\|_\infty,
\]
naturally raises the question of its relevance, namely, whether it is sharp in certain cases. In the following section, we show that for a wide variety of symbols, the new estimate provides a better bound than the classical one. In this context, we also refer the reader to Remark \ref{rem: not best}.

\section{On sharper estimates}\label{sec: estimate}

In this section, we identify a large class of self-maps $\varphi$ of $\D$ that satisfy the inequality
\[
\left\|{\frac{\sqrt{1-|\vp(0)|^2}}{1-\overline{\vp(0)}\vp}}\right\|_\infty < \sqrt{\frac{1+|\vp(0)|}{1-|\vp(0)|}},
\]
thereby establishing that the new estimate improves upon the classical one. We will eventually classify the symbols for which the above inequality holds. However, as a first observation, it follows immediately that the new estimate is sharper for symbols belonging to the open unit ball of $H^\infty(\D)$ that do not vanish at the origin:

\begin{thm}\label{thm: phi < 1}
Let $\vp$ be a self-map of $\D$. Suppose
\[
\vp(0) \neq 0,
\]
and
\[
\|\vp\|_\infty < 1.
\]
Then
\[
\|C_\vp\| \leq \left\|{\frac{\sqrt{1-|\vp(0)|^2}}{1-\overline{\vp(0)}\vp}}\right\|_\infty < \sqrt{\frac{1+|\vp(0)|}{1-|\vp(0)|}}.
\]
\end{thm}
\begin{proof}
In the proof of Theorem \ref{Bdd}, we considered the inequality in \eqref{eqn: 1 - phi bar phi}, namely
\[
|1-\overline{\vp(0)}\vp(z)| \geq 1-|\overline{\vp(0)}|,
\]
for all $z \in \D$. In the present case, in the view of $\vp(0) \neq 0$, we further observe that
\[
\|1-\overline{\vp(0)}\vp\|_\infty > 1 - |{\vp(0)}|,
\]
which completes the proof of the theorem.
\end{proof}

If a self-map $\varphi$ of $\mathbb{D}$ satisfies the condition $\|\varphi\|_\infty < 1$, as also assumed in the theorem above, then a well-known result asserts that the composition operator $C_\varphi$ is compact on $H^2(\mathbb{D})$ (see \cite[page 23]{shap}). This naturally raises the question of whether the sharp inequality obtained above continues to hold beyond the class of compact composition operators. Another question of interest concerns the case where the symbol $\varphi$ satisfies $\|\varphi\|_\infty = 1$. In Section \ref{sec: examples}, we shall observe that the scope of the sharp inequality is remarkably broad, covering all possible cases of self-maps of $\mathbb{D}$.

At present, we examine the issue of identifying self-maps of $\D$ that may result in either equality or inequality between classical and present norm estimates. We begin with simple, yet much-studied, self-maps. Fix scalars $a, b \in \mathbb{C}$ and a natural number $n \in \mathbb{N}$. Define
\[
\vp(z)=az^n+b.
\]\
Then $\vp$ is a self-map of $\D$ if and only if
\[
|b|<1 \text{ and } |a|+|b|\leq 1.
\]
Indeed, we claim that
\begin{equation}\label{eqn: phi norm a b}
\|\vp\|_\infty=|a|+|b|.
\end{equation}
This, together with the fact that $\vp(0) = b$, implies the assertion. The claim in \eqref{eqn: phi norm a b} is undoubtedly known to experts; however, since we were unable to locate it in the literature, we provide a proof here for completeness. To verify this claim, first observe, without loss of generality, that $a,b\neq 0$. Since $z^n$ maps $\D$ onto $\D$, there exists a sequence $\{z_k\} \subseteq \D$ such that
\[
\lim_{k \raro \infty} z_k^n = e^{i(\arg b-\arg a)}.
\]
This implies
\[
\begin{split}
\lim_{k \raro \infty} \varphi(z_k) & = \lim_{k \raro \infty} (az_k^n+b)
\\
& = ae^{i(\arg b-\arg a)}+b
\\
& = |a|e^{i\arg b}+b
\\
& = e^{i\arg b} (|a|+|b|).
\end{split}
\]
In particular, we have $\lim_{k \raro \infty} |\varphi(z_k)| = |a|+|b|$, and hence
\[
|a|+|b|\leq \|\varphi\|_\infty.
\]
The reverse inequality follows directly from the triangle inequality. This completes the proof of the claim. When $n=1$, this map was studied extensively in \cite[Section 3]{Deddens}. We have the following result concerning the equality between the classical and the present norm estimates:

\begin{prop}\label{AB}
Let $n\in \mathbb{N}$ and let $a, b \in \mathbb{C}$. Define $\vp(z)=az^n+b$ and assume that $a \neq 0$ and $|a|+|b|\leq 1$. Then
\[
\left\|{\frac{\sqrt{1-|\vp(0)|^2}}{1-\overline{\vp(0)}\vp(z)}}\right\|_\infty
=\sqrt{\frac{1+|\vp(0)|}{1-|\vp(0)|}},
\]
if and only if either $\vp(0)=0$ or $\|\vp\|_\infty=1$.
\end{prop}
\begin{proof}
Let $p =\left\|{\frac{\sqrt{1-|\vp(0)|^2}}{1-\overline{\vp(0)}\vp(z)}}\right\|_\infty$ and let $q = \sqrt{\frac{1+|\vp(0)|}{1-|\vp(0)|}}$. The result holds trivially if $b = 0$. Hence, without loss of generality, we assume that $b \neq 0$. Moreover, as $\vp(0) = b$, it follows that
\[
q^2 = {\frac{1-|b|^2}{(1-|b|)^2}},
\]
and also
\[p= \frac{\sqrt{1-|b|^2}}{|ab|\inf\limits_{z\in \D}  \left|z^n-\frac{1-|b|^2}{a\overline{b}}\right|}=
\frac{\sqrt{1-|b|^2}}{|ab|  \left(\frac{1-|b|^2}{|ab|}-1\right)}=
\frac{\sqrt{1-|b|^2}}{1-|b|^2-|ab|}.
\]
Therefore, we have
\[
\begin{split}
q^2 - p^2 & = (1-|b|^2)\left[ \dfrac{1}{(1-|b|)^2}-\dfrac{1}{(1-|b|^2-|ab|)^2} \right]
\\
& = \dfrac{(1-|b|^2)|b|(1-|a|-|b|)\left[ 2-|b|(1+|a|+|b|)\right]}{(1-|b|)^2(1-|b|^2-|ab|)^2} .
\end{split}
\]
Hence, $p = q$ if and only if $1-|a|-|b|=0$; that is, $\|\vp\|_\infty=1$. This completes the proof.
\end{proof}

We continue with the setting of the above proposition. Our goal is to interpret the preceding conclusion in a way that characterizes the equality in full generality. In the setting of Proposition \ref{AB}, assume that $\vp(0) \neq 0$. We claim that
\[
\|\vp\|_\infty=1,
\]
if and only if
\[
\dfrac{\vp(0)}{|\vp(0)|}\in \partial( \vp(\D)).
\]
Suppose $\|\vp\|_\infty=1$, equivalently, $|a|+|b|=1$ (see \eqref{eqn: phi norm a b}). We claim that
\[
\dfrac{b}{|b|}=\dfrac{\vp(0)}{|\vp(0)|}\in \partial( \vp(\D)).
\]
As in the proof of \eqref{eqn: phi norm a b}, there exists a sequence $\{z_k\} \subseteq \D$ such that
\[
\lim_{k \raro \infty} z_k^n = e^{i(\arg b-\arg a)}.
\]
Therefore,
\[
\begin{split}
\lim_{k \raro \infty} \vp(z_k) = ae^{i(\arg b-\arg a)}+b = |a|e^{i\arg b}+b = e^{i\arg b} (|a|+|b|) = e^{i\arg b} = \frac{b}{|b|},
\end{split}
\]
which implies
\[
\dfrac{\vp(0)}{|\vp(0)|} = \dfrac{b}{|b|} \in \partial( \vp(\D)).
\]
Conversely, assume that $\dfrac{b}{|b|}\in \partial( \vp(\D))$. If $|a|+|b|<1$, then $\vp(\overline{\D}) \subsetneqq \D$, and hence
\[
\partial( \vp(\D)) \cap \T = \emptyset,
\]
a contradiction. This completes the proof of the claim.

We now show that the above gives a complete characterization of when the classical and the new norm estimates coincide for general self-maps of $\D$.

\begin{thm}\label{equal}
Let $\vp$ be a self-map of $\D$. Then
\[
\left\|{\frac{\sqrt{1-|\vp(0)|^2}}{1-\overline{\vp(0)}\vp(z)}}\right\|_\infty = \sqrt{\frac{1+|\vp(0)|}{1-|\vp(0)|}},
\]
if and only if either
\[
\vp(0)=0,
\]
or
\[
\dfrac{\vp(0)}{|\vp(0)|}\in \partial( \vp(\D)).
\]
\end{thm}
\begin{proof}
Let $p =\left\|{\frac{\sqrt{1-|\vp(0)|^2}}{1-\overline{\vp(0)}\vp(z)}}\right\|_\infty$ and let $q = \sqrt{\frac{1+|\vp(0)|}{1-|\vp(0)|}}$. If $\varphi(0) = 0$, then the equality $p = q$ holds trivially. Assume now that $\varphi(0) \neq 0$. The equality $p=q$ can be rewritten as
\[
\frac{1}{1-|\vp(0)|}= \left\|{\frac{1}{1-\overline{\vp(0)}\vp(z)}}\right\|_\infty.
\]
As
\[
\left\|{\frac{1}{1-\overline{\vp(0)}\vp(z)}}\right\|_\infty  = \frac{1}{\inf\limits_{z\in \D}|1-\overline{\vp(0)}\vp(z)|},
\]
the equality $p=q$ is equivalent to
\[
1-|\vp(0)|=\inf\limits_{z\in \D}|1-\overline{\vp(0)}\vp(z)|.
\]
This happens if and only if there exists a sequence $\{z_n\} \subseteq \D$ such that
\begin{equation}\label{l}
\lim_{n\rightarrow \infty}|1-\overline{\vp(0)}\vp(z_n)|=1-|\vp(0)|.
\end{equation}
By the Bolzano–Weierstrass theorem, there exists a subsequence of $\{z_n\}$, denoted again by $\{z_n\}$, such that
$\{\vp(z_n)\}$ converges to some number, say $l$. By \eqref{l}, we have
\[
|1-\overline{\vp(0)}l|=1-|\vp(0)|.
\]
Using the triangle inequality, we find
\[
1-|\overline{\vp(0)}l|\leq 1-|\vp(0)|.
\]
Since $\varphi(0) \neq 0$, this implies
\[
|\overline{\vp(0)}l|\geq |\vp(0)|,
\]
and hence
\[
|l|\geq 1.
\]
On the other hand, as $|\vp(z_n)| \leq 1$ for all $n$, it follows that $|l| \leq 1$. Therefore, we conclude that
\[
|l|=1.
\]
On the other hand, by squaring both sides of $|1-\overline{\vp(0)}l|=1-|\vp(0)l|$, we see that
\[
\overline{\vp(0)}l=|\vp(0)|.
\]
Hence,
\[
\dfrac{\vp(0)}{|\vp(0)|}=l\in \partial( \vp(\D)).
\]
For the converse direction, if we assume that $l\in \partial( \vp(\D))$, then the above computation can be reversed, leading back to \eqref{l}, and this completes the proof.
\end{proof}

In Section \ref{sec: examples}, we illustrate the new norm estimate through a few specific examples.

\section{Reproducing kernel Hilbert spaces}\label{sec-RKHS}

Following the ideas used in the proof of Theorem \ref{Bdd}, we now establish a general result on the boundedness of composition operators, along with a corresponding norm estimate, in the setting of reproducing kernel Hilbert spaces.

We consider kernel functions $k_1$ and $k_2$ defined on two sets $X_1$ and $X_2$, respectively. Let $\vp: X_2\rightarrow X_1$ be a function. The composition operator $C_\vp$ on $\clh(k_1)$ is defined as usual by
\[
C_\vp f = f \circ \vp,
\]
for all $f \in \clh(k_1)$. We recall a standard fact about composition operators on reproducing kernel Hilbert spaces (cf. \cite[ page 72, Theorem 5.10]{paul}, which will be used in the general result below: the composition operator $C_\vp$ from $\clh(k_1)$ to $\clh(k_2)$ is bounded if and only if there exists a constant $c > 0$ such that
\begin{equation}\label{eqn: gen CO norm}
c^2 k_1(x,y) - k_2(\vp(x), \vp(y)) \geq 0.
\end{equation}
Moreover, in this case, $\|C_\vp\|$ is the least such constant $c$.

\begin{thm}\label{genbdd}
Let $k_1$ and $k_2$ be kernel functions on $X_1$ and $X_2$, respectively, and let $\vp: X_2\rightarrow X_1$ be a function. Define $k :X_2\times X_2\rightarrow\C$ by
\[
k(x,y)= \frac{k_2(x,y)}{k_1(\vp(x),\vp(y))}.
\]
Assume that $k$ is a kernel function on $X_2$. Suppose further that for some fixed $\alpha \in X_2$,
\[
\frac{1}{k(\cdot, \alpha)} \in \clm(\clh(k_2), \clh(k_2)).
\]
Then $C_\vp$ is a bounded linear operator from $\clh(k_1)$ to $\clh(k_2)$. Moreover,
\[
\|C_\vp\|\leq \sqrt{k(\alpha,\alpha)}\left\|M_{\frac{1}{k(\cdot, \alpha)}}\right\|.
\]
\end{thm}
\bpf
We follow the lines of argument of the proof of Theorem \ref{Bdd}. Since $k(\cdot, \alpha) \in \clh(k)$, applying \eqref{eqn: f in RKHS}, we conclude
\[
\displaystyle\| k(\cdot, \alpha) \|_{\clh(k)}^2 k(x,y) - k(x, \alpha) \overline{k(y, \alpha)} \geq 0.
\]
Multiplying this with kernel functions $k_1(\vp(x),\vp(y))$ and $\displaystyle\frac{1}{k(x, \alpha) \overline{k(y, \alpha)}}$, we obtain
\begin{equation}\label{E2}
\displaystyle\| k(\cdot, \alpha)\|_{\clh(k)}^2\frac{k_2(x,y)}{k(x, \alpha) \overline{k(y, \alpha)}} - k_1(\vp(x),\vp(y))\geq 0.
\end{equation}
Given that $\frac{1}{k(\cdot, \alpha)}$ is a multiplier in $\clm(\clh(k_2), \clh(k_2))$, \eqref{eqn: mult} implies
\begin{equation}\label{E3}
\displaystyle\left\|M_{\frac{1}{k(\cdot, \alpha)}}\right\|^2k_2(x,y)-\frac{k_2(x,y)}{k(x, \alpha) \overline{k(y, \alpha)}}\geq0.
\end{equation}
Finally, multiplying the previous inequality by the positive constant $\| k(\cdot, \alpha)\|_{\clh(k)}^2$ and adding it to \eqref{E2} yields
\[
\| k(\cdot, \alpha)\|_{\clh(k)}^2 \left\|M_{\frac{1}{k(\cdot, \alpha)}} \right\|^2k_2(x,y)-k_1(\vp(y),\vp(x))\geq 0.
\]
This, in view of \eqref{eqn: gen CO norm}, implies that $C_\vp$ is a bounded linear operator and
\[
\|C_\vp\| \leq \| k(\cdot, \alpha)\|_{\clh(k)} \left\|M_{\frac{1}{k(\cdot, \alpha)}} \right\| = \sqrt{k(\alpha,\alpha)}\left\|M_{\frac{1}{k(\cdot, \alpha)}}\right\|,
\]
which completes the proof of the theorem.
\epf

Composition operators between reproducing kernel Hilbert spaces are interesting objects. The above result can be applied, for instance, to prove the boundedness of composition operators between $H^2(\D)$ and the Hardy space $H^2(\B^n)$ over the unit ball $\B^n \subseteq \mathbb{C}^n$.

\begin{cor}
Let $\vp:\B^n\rightarrow\D$ be an analytic function. Then $C_\vp$ is a bounded linear operator from $H^2(\D)$ to $H^2(\B^n)$.
\end{cor}
\bpf
Following the setting of Theorem \ref{genbdd}, we write the kernel for $H^2(\D)$ as $k_1(z,w)=\frac{1}{1-z \overline{w}}$, and that for $H^2(\B^n)$ as
\[
k_2(\boldsymbol{z}, \boldsymbol{w})=\displaystyle \frac{1}{(1-\langle \boldsymbol{z}, \boldsymbol{w}  \rangle)^n},
\]
where $\langle \boldsymbol{z}, \boldsymbol{w}  \rangle = \sum_{i=1}^n z_i \bar{w}_i$ for all $\boldsymbol{z}, \boldsymbol{w} \in \B^n$. Also, choose
\[
\alpha = \boldsymbol{0} \in \mathbb{C}^n.
\]
In this case, we have
\[
k(\boldsymbol{z}, \boldsymbol{w})=\frac{k_2(\boldsymbol{z}, \boldsymbol{w})}{k_1(\vp(\boldsymbol{z}),\vp(\boldsymbol{w}))}= \frac{1- \vp(\boldsymbol{z}) \overline{\vp(\boldsymbol{w})}}{(1-\langle \boldsymbol{z}, \boldsymbol{w} \rangle)^n}.
\]
Clearly, $k(\boldsymbol{z}, \boldsymbol{w}) \geq 0$ on $\B^n$. Moreover, it is easy to see that
\[
\frac{1}{k(\cdot, \boldsymbol{0})}(\boldsymbol{z})=\frac{1}{1-\overline{\vp(\boldsymbol{0})}\vp(\boldsymbol{z})},
\]
is a bounded analytic function on $\B^n$, and hence, a multiplier of $H^2(\B^n)$. Therefore, by Theorem \ref{genbdd}, the composition operator $C_\vp$ maps $H^2(\D)$ boundedly into $H^2(\B^n)$.
\epf

Along similar lines, one can also show that for any analytic function $\vp: \D_n\rightarrow\D$, the composition operator $C_\vp:H^2(\D)\rightarrow H^2(\D^n)$ is always bounded, where $H^2(\D^n)$ denotes the Hardy space over the polydisc $\D^n$. In this context, we refer the reader to \cite[Proposition 3]{jaf}.

\section{Examples}\label{sec: examples}

As pointed out in Theorem \ref{thm: phi < 1}, if $\varphi$ is a self-map of $\mathbb{D}$ and $\|\varphi\|_\infty < 1$, then the new norm estimate for $C_\varphi$ is sharper than the classical one. It is now a natural question of interest whether self-maps with
\[
\|\varphi\|_\infty = 1,
\]
can also yield such a sharper estimate. In the following example, we show that there are, in fact, many such cases.

\begin{ex}\label{eg}
\begin{figure}[h]
        \includegraphics[width=0.5\textwidth]{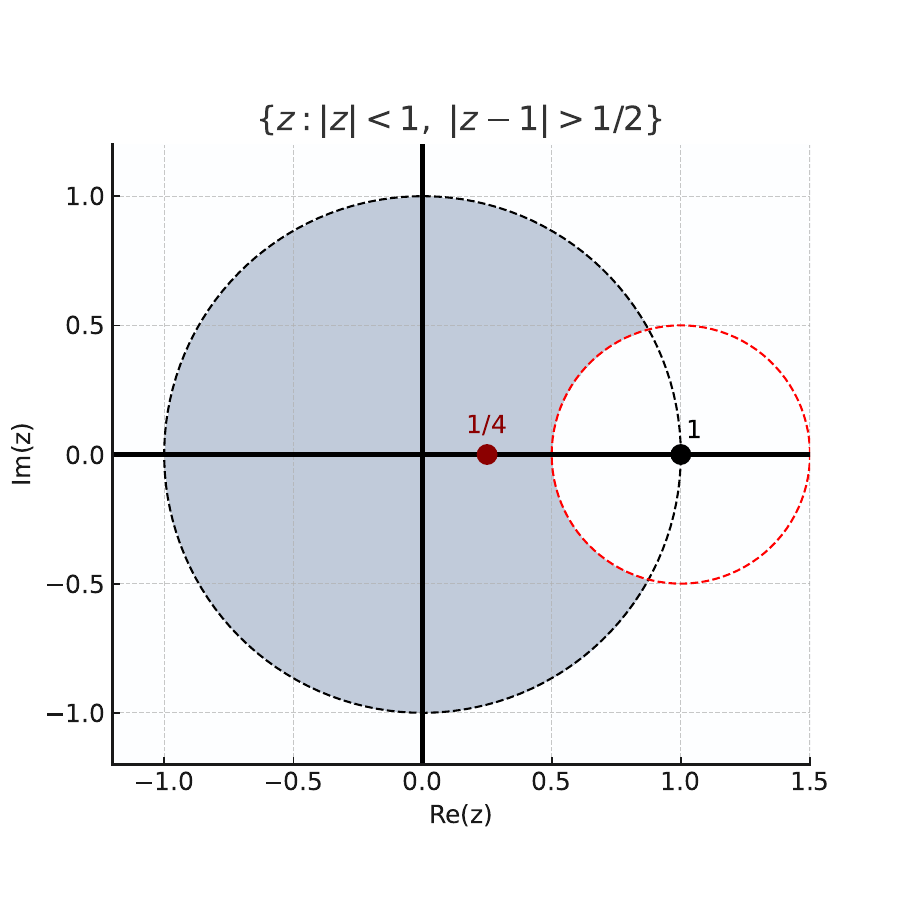}
    \end{figure}
Let
\[
\Omega=\{z\in \D: |z-1|>1/2 \}.
\]

By the Riemann mapping theorem, we can choose a biholomorphic map $\vp$ from $\D$ onto $\Omega$ with
\[
\vp(0)=\frac{1}{4}.
\]
Clearly, $\|\vp\|_\infty=1$ by definition of $\vp$. Since $\vp(0)=\frac{1}{4} > 0$, we have
\[
\dfrac{\vp(0)}{|\vp(0)|}=1,
\]
But
\[
\overline{\vp(\D)} = \overline{\Omega},
\]
and $1 \notin \overline{\Omega}$. This, in particular, implies that $1 \notin \partial \overline{\Omega}$. Hence, by Theorem \ref{equal}, we conclude that
\[
\left\|{\frac{\sqrt{1-|\vp(0)|^2}}{1-\overline{\vp(0)}\vp(z)}}\right\|_\infty <\sqrt{\frac{1+|\vp(0)|}{1-|\vp(0)|}} \cdot
\]
\end{ex}

The above example is particularly noteworthy: although the function $\varphi$ is not inner, it satisfies
\[
|\vp(e^{it})| = 1,
\]
on a subset of $\mathbb{T}$ of positive Lebesgue measure. This, in particular, implies that $C_\vp$ is a noncompact operator on $H^2(\mathbb{D})$. This feature is especially relevant in view of the fact that $C_\varphi$ is always compact whenever
\[
\|\vp\|_\infty < 1.
\]

In this context, the following example is also noteworthy, as its symbol has norm one while the corresponding composition operator is compact.
\begin{figure}[h]
        \includegraphics[width=0.5\textwidth]{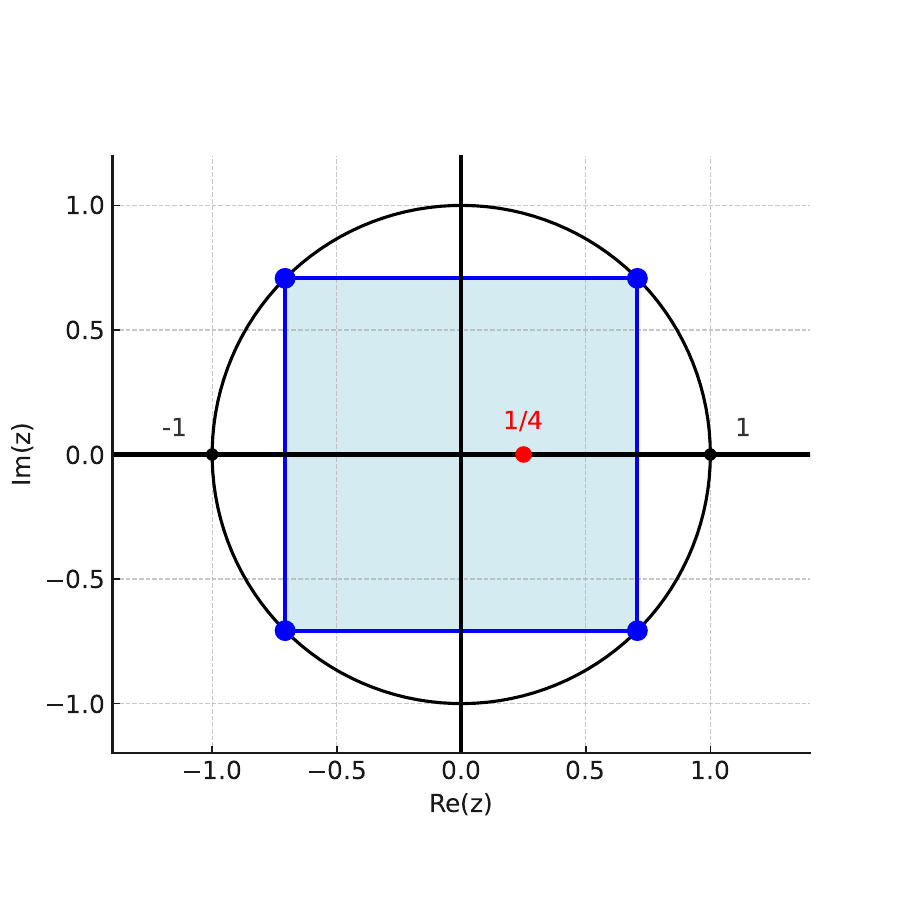}
    \end{figure}
\begin{ex}\label{examp: 1 and i}
Let $\Omega$  be the interior of the square inscribed in the unit disk with vertices at $\frac{1}{\sqrt{2}}(\pm 1 \pm i)$.
Again, by the Riemann mapping theorem, we can choose a biholomorphic map $\vp$ from $\D$ onto $\Omega$ such that
\[
\vp(0)=\frac{1}{4}.
\]
Clearly, by the definition of $\vp$, we have $\|\vp\|_\infty=1$. Furthermore, by the polygonal compactness theorem \cite[Chapter 2]{shap},
the operator $C_\vp$ is compact on $H^2(\D)$.
As in the previous example, since $\vp(0)=\frac{1}{4}$, we have $\dfrac{\vp(0)}{|\vp(0)|}=1$, which is not contained in $\overline{\vp(\D)}$. Hence, by Theorem \ref{equal}, it follows that
\[
\left\|{\frac{\sqrt{1-|\vp(0)|^2}}{1-\overline{\vp(0)}\vp(z)}}\right\|_\infty
<\sqrt{\frac{1+|\vp(0)|}{1-|\vp(0)|}} \cdot
\]
In other words, there exist a self-map $\vp$ of $\D$ such that $C_\vp$ is compact on $H^2(\D)$, $\|\vp\|_\infty=1$, and the new estimate is sharper than the classical one.
\end{ex}

The following example is particularly concrete, involving the computation for a degree-one polynomial self-map of the unit disc. This example is classical and has been studied extensively by Deddens \cite{Deddens}. The idea here is not to apply Theorem \ref{equal} directly, but rather to proceed by direct computation.

\begin{ex}\label{aff}
Define
\[
\vp(z)= \frac{1}{3} (1+z).
\]
Clearly, $\vp$ is a self-map of $\D$. Moreover, as $|\vp(z)|\leq \frac{2}{3}$ on $\D$, it follows that
\[
|1-\overline{\vp(0)}\vp(z)| \geq 1- |\overline{\vp(0)}|\ |\vp(z)|\geq 1-\frac{1}{3}\frac{2}{3}=\frac{7}{9},
\]
for all $z\in \D$. Therefore,
\[
\left\|{\frac{\sqrt{1-|\vp(0)|^2}}{1-\overline{\vp(0)}\vp(z)}}\right\|_\infty\leq \frac{9}{7}\sqrt{\frac{8}{9}} =\frac{6\sqrt{2}}{7},
\]
and hence
\[
\|C_{\vp}\|\leq \left\|{\frac{\sqrt{1-|\vp(0)|^2}}{1-\overline{\vp(0)}\vp(z)}}\right\|_\infty
\leq\frac{6\sqrt{2}}{7}\approx 1.212<1.414\approx\sqrt{2}=\sqrt{\frac{1+|\vp(0)|}{1-|\vp(0)|}}.
\]
This shows that the new estimate is indeed sharper than the classical one.
\end{ex}

We proceed to further elaborate on the exotic nature of composition operators in the context of our results.

\begin{ex}
We continue with the domain $\Omega$ introduced in Example \ref{examp: 1 and i}. This time, we choose a biholomorphic map $\vp$ from $\D$ onto $\Omega$ such that
\[
\vp(0)= \frac{1}{4}+i\frac{1}{4}.
\]
By the polygonal compactness theorem, $C_{\vp}$ is compact. Moreover,
\[
\dfrac{\vp(0)}{|\vp(0)|} = \dfrac{1}{\sqrt{2}}(1+i)\in \partial \overline{\Omega},
\]
and hence Theorem \ref{equal} implies that
\[
\left\|{\frac{\sqrt{1-|\vp(0)|^2}} {1-\overline{\vp(0)}\vp(z)}}\right\|_\infty =\sqrt{\frac{1+|\vp(0)|}{1-|\vp(0)|}} \cdot
\]
\end{ex}

Finally, we present an example of a non-compact composition operator for which the classical and the new norm estimates coincide:

\begin{ex}
Along the lines of Deddens \cite{Deddens} and Example \ref{aff}, we consider the self-map of $\D$ as
\[
\displaystyle\vp(z)=\frac{1+z}{2}.
\]
Then
\[
\dfrac{\vp(0)}{|\vp(0)|}=1\in\overline{\vp(\D)},
\]
and thus Theorem \ref{equal} implies that
\[
\left\|{\frac{\sqrt{1-|\vp(0)|^2}}{1-\overline{\vp(0)}\vp(z)}}\right\|_\infty =\sqrt{\frac{1+|\vp(0)|}{1-|\vp(0)|}} \cdot
\]
However, $C_\vp$ is not compact \cite[Page 30]{shap}.
\end{ex}

In summary, for all self-maps $\vp$ with $\|\vp\|_\infty<1$ and $\vp(0)\neq 0$, we know that $C_\vp$ is compact, and
\[
\left\|{\frac{\sqrt{1-|\vp(0)|^2}}{1-\overline{\vp(0)}\vp(z)}}\right\|_\infty < \sqrt{\frac{1+|\vp(0)|}{1-|\vp(0)|}}.
\]
In contrast, when $\|\varphi\|_\infty = 1$, a definitive general conclusion is elusive: all possible combinations, namely, $C_\vp$ being either compact or non-compact, and the inequality between the classical and the present norm estimates being either strict or an equality, may occur.

We conclude this section with the following remark concerning the new norm estimate obtained in this paper.

\begin{rem}\label{rem: not best}
Let $\vp$ be a self-map of $\D$. Recall that \cite[Corollary 3.7]{cow}
\[
\sqrt{\frac{1}{1-|\vp(0)|^2}}\leq \|C_\vp\| \leq \sqrt{\frac{1+|\vp(0)|}{1-|\vp(0)|}}.
\]
It is easy to verify that the lower bound is attained by composition operators induced by constant self-maps, while it is known that the upper bound is attained by those induced by inner functions \cite{Nord}. Therefore, if $\vp$ is inner, then
\[
\|C_\vp\|=\left\|{\frac{\sqrt{1-|\vp(0)|^2}}{1-\overline{\vp(0)}\vp(z)}}\right\|_\infty=\sqrt{\frac{1+|\vp(0)|}{1-|\vp(0)|}}.
\]
On the other hand, if $\vp\equiv c$ is a constant function, then
\[
\sqrt{\frac{1}{1-|c|^2}}= \|C_\vp\|=\left\|{\frac{\sqrt{1-|\vp(0)|^2}}{1-\overline{\vp(0)}\vp(z)}}\right\|_\infty.
\]
Thus, our estimate agrees with $\|C_\vp\|$ in both extreme cases. This naturally leads to the question of whether, for an arbitrary self-map $\vp$ of $\D$, the norm of $C_\vp$ always coincides with the new estimate; that is, whether
\[
\|C_\vp\|=\left\|{\frac{\sqrt{1-|\vp(0)|^2}}{1-\overline{\vp(0)}\vp(z)}}\right\|_\infty,
\]
holds for all self-maps $\vp$ of $\D$. However, this is not true in general. For instance, consider the self-map of $\D$ discussed in Example \ref{aff}. In that case, the exact value of $C_\vp$ is known explicitly \cite[Theorem 9.4]{cow}:
\[
\|C_{\vp}\|=\sqrt{\frac{6}{3+\sqrt{5}}}\approx 1.07.
\]
This ensures that the new estimate provides a genuinely general upper bound and need not be norm always.
\end{rem}

\section{H. Schwartz’s estimates}\label{sec: Schwartz}

In Section 3 of his thesis \cite{Sch}, Schwartz claimed an improved norm estimate for composition operators compared to the classical one given by Littlewood’s theorem (see the introduction to Section 3, page 30 of \cite{Sch} and then \cite[Theorem 3.10]{Sch}). His estimate also involves the supremum norm of functions. More specifically (the $p=2$ version of \cite[Theorem 3.10]{Sch}): Let $\vp$ be a self-map of $\D$. Assume that $\|\vp\|_\infty < 1$ and $\vp(0) \neq 0$. Then
\[
\|C_\vp\| < \sqrt{\frac{\|\vp\|_\infty + |\vp(0)|}{\|\vp\|_\infty - |\vp(0)|}}.
\]
This naturally leads to a question regarding its relation to the new norm estimate established in this paper. We address this issue as follows:

\begin{prop}\label{prop: Schartz}
Let $\vp$ be a self-map of $\D$. If $\|\vp\|_\infty<1$ and $\vp(0)\neq0$, then
\[
\dfrac{1+|\vp(0)|}{1-|\vp(0)|} < \dfrac{\|\vp\|_\infty+|\vp(0)|}{\|\vp\|_\infty-|\vp(0)|}.
\]
\end{prop}
\bpf
Set $\|\vp\|_\infty=r$ and $|\vp(0)|=a$. Then a straightforward calculation shows that
\[
\dfrac{1+a}{1-a}<\dfrac{r+a}{r-a},
\]
if and only if
\[
a(1-r)> 0,
\]
which indeed holds since $0<a<r<1$.
\epf

Therefore, for a self-map $\vp$ of $\D$, if $\|\vp\|_\infty<1$ and $\vp(0)\neq0$, then
\[
\|C_{\vp}\|\leq \left\|{\frac{\sqrt{1-|\vp(0)|^2}}{1-\overline{\vp(0)} \vp}}\right\|_\infty < \sqrt{\frac{1+|\vp(0)|}{1-|\vp(0)|}} < \sqrt{\dfrac{\|\vp\|_\infty+|\vp(0)|}{\|\vp\|_\infty-|\vp(0)|}}.
\]
This means that not only is the new norm estimate sharper, but even the classical norm estimate in Littlewood’s theorem turns out to be sharper than Schwartz’s estimate. This is clearly contrary to the claim made by Schwartz \cite[Section 3, page 30]{Sch}.

As a concluding remark, we note that composition operators constitute one of the most classical families of operators in functional analysis and are connected with a wide variety of subjects. As previously pointed out, composition operators are of particular interest when studied on different function spaces, such as Besov spaces \cite{ero}, weighted Hardy spaces \cite{Herve}, and spaces in several variables \cite{kos}. They are also closely related to other areas, including semigroup theory and the invariant subspace problem \cite{eva}. Further investigation is required to determine whether reproducing kernel methods can meaningfully supplement these theories. It would also be worthwhile to study whether the reproducing kernel methods developed here can be used to compute the essential norm of composition operators (see \cite{Joel}), particularly in several variable settings.

\vspace{0.2in}

\noindent\textbf{Acknowledgement:}
The first author acknowledges the support of the Prime Minister's Research Fellowship scheme (PMRF ID. 2303407). The research of the third named author is supported in part by TARE (TAR/2022/000063) by SERB, Department of Science \& Technology (DST), Government of India.

\end{document}